\documentclass[12pt]{article}

\newcommand{\captionfonts}{\footnotesize} 

			\makeatletter  
			\long\def\@makecaption#1#2{%
			  \vskip\abovecaptionskip
			  \sbox\@tempboxa{{\captionfonts #1: #2}}%
			  \ifdim \wd\@tempboxa >\hsize
			    {\captionfonts #1: #2\par}
			  \else
			    \hbox to\hsize{\hfil\box\@tempboxa\hfil}%
			  \fi
			  \vskip\belowcaptionskip}
			\makeatother   

\usepackage{amsmath, amsthm, amssymb}
\usepackage{amsfonts}
\usepackage{epsfig}
\usepackage{graphics}
\usepackage{subfigure}
\usepackage{graphicx}
\usepackage{color}
\usepackage{epstopdf} 
\usepackage{float} 
\usepackage{mathrsfs} 
\usepackage{relsize} 
\newcommand{\ds}{\displaystyle}

\newcommand{\NN}{\mathbb{N}}

\newtheorem{thm}{Theorem}
\newtheorem{df}{Definition}
\newtheorem{rem}{Remark}
\newcommand{\bq}{\begin{equation}}
\newcommand{\eq}{\end{equation}}
\newcommand{\bqr}{\begin{eqnarray}}
\newcommand{\eqr}{\end{eqnarray}}
\newcommand{\bqrn}{\begin{eqnarray*}}
\newcommand{\eqrn}{\end{eqnarray*}}

\begin{document}
\begin{center}	
{\bf Substitution Method for Fractional Differential Equations}\\[2mm] 
Pavel B. Dubovski, Jeffrey A. Slepoi
\end{center}

\pagenumbering{arabic}
\setcounter{page}{1} 

\thispagestyle{empty}


\pagenumbering{roman}
\setcounter{page}{1} 

\thispagestyle{plain}
\pagenumbering{arabic}

\begin{abstract}
	Numerical solving differential equations with fractional derivatives requires elimination of the singularity which is inherent in the standard definition of fractional derivatives.  The method of integration by parts to eliminate this singularity is well known. It allows to solve some equations but increases the order of the equation and sometimes leads to wrong numerical results or instability.
	We suggest another approach: the elimination of singularity by substitution. It does not increase the order of equation and its numerical implementation provides the opportunity to define fractional derivative as the limit of discretization. 
		We present a sufficient condition for the substitution-generated difference approximation to be well-conditioned.
	We demonstrate how some equations can be solved using this method with full confidence that the solution is accurate with at least second order of approximation.         		
\end{abstract}
 {\it MSC 2010\/}: 26A33, 65L20\\
Keywords: fractional differential equations; substitution method; well-conditioned problem; numerical stability; numerical accuracy; approximation

\section{Introduction}
There are multiple definitions of fractional derivatives, of which the best known are the Riemann-Liouville and Caputo definitions.  Caputo's definition of fractional derivative of order $\alpha>0$ is especially well-suited for initial value problems, which we mostly consider in this paper:
\begin{equation}\label{CaputoFD}
D^\alpha f(t) := \frac{1}{\Gamma(n-\alpha)}\int_{a}^{t} \frac{f^{(n)}(x)}{(t-x)^{\alpha-n+1}}dx, 
\end{equation}
where $\Gamma$ is the gamma function and $n-1 \le \alpha < n, n \in \mathbb{N}$. Thus, $n=\lceil \alpha \rceil$.
It is worth pointing out that Caputo derivative was first derived in 1948 by Gerasimov \cite{Gerasimov} and, independently, by Caputo in 1967 \cite{Caputo}. 

As an example of another definition, which possesses certain good properties, we can mention the definition of Khalil et al \cite{2014Khalil}. However, the nature of this interesting definition is local and, like the usual integer derivatives, depends only on "present", whereas Caputo and Riemann-Liouville derivatives are essentially nonlocal and, thus, provide additional tools for modeling the systems, depending on 
"past". Among such systems we can mention diffusion, chaos and anomalous transport \cite{Edelman, Zaslavsky}. More information about physical backgrounds for differential calculus can be found, e.g., in \cite{Edelman, Hilfer}.\\

This paper deals with differential equations containing Caputo and Riemann-Liouville derivatives. 
As we see in (\ref{CaputoFD}), the main computational obstacle is the presence of mild singularity.  Integrating the fractional derivative by parts helps to fix this problem.  However, in some cases, the use of finite difference method to the equation modified by integration by parts does not guarantee convergence to the correct solution 
\cite{2020-1Du-Sl}. 

In this paper we present another approach to eliminate the singularity: the method of substitution. We show that it provides an opportunity to view
fractional derivative differently and leads us to another way to solve fractional differential equations.    

Similar discretization of Caputo fractional derivative can be achieved by splitting the derivative into multiple intervals, assuming that the regular derivative is constant under the integral and integrating the remaining part as was shown in \cite{Sousa} for $1<\alpha<2$.  We present a more general approach, which is possible through substitution and prove the theorem that substitution method implemented by the trapezoidal numerical integration, guarantees the order of approximation greater than 1 for any order of fractional derivatives $\alpha > 0$.   

On the basis of the substitution method, we present the difference schemes of the second order of approximation for linear equations for any $\alpha>0$:  section \ref{SectionDiscr1} for $0<\alpha<1$, section \ref{SectionDiscr2} for $0<\alpha<2$, and section \ref{SectionDiscrN} for any $\alpha >0$.

In section \ref{SectionSuffCond} we provide the sufficient condition for the corresponding (substitution-based)  finite difference scheme to be well-conditioned.

Finally, in section \ref{SectionExamples}, we present a few examples of calculations for (a) an equation, previously considered in \cite{Podlubny} using the Gr\"{u}nwald-Letnikov derivative, and (b) for the hyper-generalized fractional Bessel equation \cite{2020-4Du-Sl}, for which we know the exact solution. These numerical  results demonstrate the strength and exactness of the proposed substitution method.

\section {Elimination of singularity in fractional derivative}\label{SectionSubstMeth}
In this section we eliminate the singularity in the Caputo fractional derivative and prove the validity of its alternative formula.

We would like to eliminate the singularity in \eqref{CaputoFD}.  We introduce variable $u$ such that 
\begin{equation*}
dx = K\cdot(t-x)^{\alpha-n+1}du,
\end{equation*}
where $K$ is a constant.  Then
\begin{equation*}
u = \frac{(t-x)^{n-\alpha}}{K\cdot(n-\alpha)}. 
\end{equation*}
For simplification, let $K=\frac{1}{n-\alpha}$, then
\begin{equation*}
u = (t-x)^{n-\alpha} \text { and } x =t - u^{\frac{1}{n-\alpha}}. 
\end{equation*}
This yields:
\begin{eqnarray}
du &=& -(n-\alpha)(t-x)^{n-\alpha-1}dx, \nonumber \\
dx &=& -\frac{1}{n-\alpha}\cdot (t-x)^{\alpha-n+1} du, \nonumber \\
f^{(n)}(x)&=& f^{(n)}(t-u^{\frac{1}{n-\alpha}}). \nonumber 
\end{eqnarray}
Therefore
\begin{eqnarray}\label{NoSingForm}
D^\alpha f(t)&  = &\frac{1}{\Gamma(n-\alpha)}\int_{0}^{t} \frac{f^{(n)}(x)}{(t-x)^{\alpha-n+1}}dx\nonumber \\
 & = &-\frac{1}{\Gamma(n-\alpha)}\cdot\frac{1}{n-\alpha} \int_{t^{n-\alpha}}^{0}f^{(n)}(t-u^{\frac{1}{n-\alpha}})du \nonumber \\
 & = &\frac{1}{\Gamma(n+1-\alpha)} \int_{0}^{t^{n-\alpha}}f^{(n)}(t-u^{\frac{1}{n-\alpha}})du. \nonumber 
\end{eqnarray}
These results provide another representation of the Caputo fractional derivative:
\begin{equation}\label{SubstDefFrD}
D^\alpha f(t) = \frac{1}{\Gamma(n+1-\alpha)} \int_{0}^{t^{n-\alpha}}f^{(n)}(t-u^{\frac{1}{n-\alpha}})du.
\end{equation}

The obtained representation of the derivative can be numerically calculated using the trapezoidal method by calculating $f^{(n)}(x_k)$ at each point $x_k=x_{k-1}+h_k$ between $0$ and $t$ ($h_k$ is a step from $x_{k-1}$ to $x_k$) and recalculating $\Delta_k=(t-x_{k-1})^{n-\alpha}-(t-x_k)^{n-\alpha}>0$. The step $h_k$ between the points $x_{k-1}$ and $x_k, k=0,1,...,M$ does not have to be constant, but if we have several derivatives in equation, there should be only one grid for all of them, which is based on the grid spacing $h_k$ along $x$. Based on the trapezoidal rule, we obtain from \eqref{SubstDefFrD}:
\begin{equation}\label{FrDSubst}
D^\alpha f(t) \approx \frac{1}{\Gamma(n+1-\alpha)} \sum_{k=1}^{m}\frac{f^{(n)}(x_k)+f^{(n)}(x_{k-1})}{2}\cdot(u_{k-1}-u_{k}).
\end{equation}
The finite difference approximation of each derivative is based on $u_k=(t-x_k)^{n-\alpha}$.  Since the increase in $x$ leads to the decrease in $u$, $u_{k-1}>u_k$.  We used approximation of the integral in \eqref{SubstDefFrD} by trapezoidal method because it is a simple representation and hence provides a straight forward numerical implementation of Caputo derivatives. 
	\begin{rem}
	Riemann-Liouville fractional derivative 
	\begin{equation*}
	D^\alpha_{RL} f(t) = \frac{1}{\Gamma(n-\alpha)}\frac{d^n}{dt^n}\int_{0}^{t}\frac{f(x)}{(t-x)^{\alpha+1-n}}dx, n-1\le\alpha<n,
	\end{equation*}
	through identical steps can be represented as
\begin{equation*}
D^\alpha_{RL} f(t) = \frac{1}{\Gamma(n+1-\alpha)}\frac{d^n}{dt^n}\int_{0}^{t^{n-\alpha}}f(t-u^{\frac{1}{n-\alpha}})du.
\end{equation*} 
If we take the derivative, we get the well-known formula \cite{Podlubny} which connects Riemann-Liouville and Caputo derivatives for $f(t)\in C^n$. For example, for $0<\alpha<1$ we have
\begin{eqnarray}
D^\alpha_{RL} f(t) &=& \frac{1}{\Gamma(2-\alpha)}\frac{d}{dt}\int_{0}^{t^{1-\alpha}}f(t-u^{\frac{1}{1-\alpha}})du\nonumber \\
&=&	\frac{1}{\Gamma(2-\alpha)}f(0)(1-\alpha)t^{-\alpha}+\frac{1}{\Gamma(2-\alpha)} \int_{0}^{t^{1-\alpha}}f'(t-u^{\frac{1}{1-\alpha}})du =\frac{f(0)t^{-\alpha}}{\Gamma(1-\alpha)}+D^\alpha f(t). \nonumber 
\end{eqnarray}
\end{rem} 
As we see, this remark indicates that the proposed substitution method also can be 
successfully used for solving fractional differential equations in the Riemann-Liouville sense.

%
\begin{thm}{}\label{Thm1}
The Caputo derivative formula \eqref{SubstDefFrD} is equivalent to
\begin{equation}
D^\alpha f(t) =\frac{1}{\Gamma(n+1-\alpha)}\lim_{h_{\max}\to 0} \sum_{k=1}^{m}\frac{f^{(n)}(x_k)+f^{(n)}(x_{k-1})}{2}\cdot\left((t-x_{k-1})^{n-\alpha}-(t-x_k)^{n-\alpha}\right), \label{4}
\end{equation}
where $0=x_0<x_1<...<x_m=t, h_k=x_k-x_{k-1}>0, h_{\max}=\max\limits_{k}\{h_k\},  h_{\max}\underset{m\to \infty}{\longrightarrow} 0$ and \\ $m\cdot h_{\max} \le C(t)<\infty$.  Moreover, order of approximation of Caputo derivative through \eqref{FrDSubst} is at least $O(h_{\max}^{n-\alpha+1})$.
\end{thm}
\begin{proof}
By using the grid and assuming 
\begin{equation*}
D_k^\alpha f(t) =\int_{(t-x_k)^{n-\alpha}}^{(t-x_{k-1})^{n-\alpha}}f^{(n)}(t-u^{\frac{1}{n-\alpha}})du,
\end{equation*}
we obtain from \eqref{SubstDefFrD}
	\begin{equation*}
	D^\alpha f(t) = \frac{1}{\Gamma(n+1-\alpha)}\sum_{k=1}^{m}D_k^\alpha f(t).
	\end{equation*} 
Now, let's consider the residual $R$, which is the difference of \eqref{SubstDefFrD} and \eqref{FrDSubst}:
\begin{subequations}
	\begin{align*}	
 R&=\frac{1}{\Gamma(n+1-\alpha)} \int_{0}^{t^{n-\alpha}}f^{(n)}(t-u^{\frac{1}{n-\alpha}})du \\ 
&\phantom{-} -\frac{1}{\Gamma(n+1-\alpha)}\sum_{k=1}^{m}\frac{f^{(n)}(x_k)+f^{(n)}(x_{k-1})}{2}\left((t-x_{k-1})^{n-\alpha}-(t-x_k)^{n-\alpha}\right) \\
 &=  \frac{1}{\Gamma(n+1-\alpha)}\sum_{k=1}^{m}\Big[\int_{(t-x_k)^{n-\alpha}}^{(t-x_{k-1})^{n-\alpha}}f^{(n)}(t-u^{\frac{1}{n-\alpha}}) du \\
&\phantom{-} -\frac{f^{(n)}(x_k)+f^{(n)}(x_{k-1})}{2}\left((t-x_{k-1})^{n-\alpha}-(t-x_k)^{n-\alpha}\right)\Big].
	\end{align*}
\end{subequations}
From the mean value theorem we know that
\begin{equation}\label{MeanValThrm}
D^\alpha_k f(t)=\int_{(t-x_k)^{n-\alpha}}^{(t-x_{k-1})^{n-\alpha}}f^{(n)}(t-u^{\frac{1}{n-\alpha}}) du = f^{(n)}(t-u_j^{\frac{1}{n-\alpha}})\left((t-x_{k-1})^{n-\alpha}-(t-x_k)^{n-\alpha}\right),
\end{equation}
where $(t-x_k)^{n-\alpha}<u_j<(t-x_{k-1})^{n-\alpha}$.  Relying upon $u_j$, we choose the corresponding value of $x_j$ from the expression $u_j=(t-x_j)^{n-\alpha}$.  Then, we obtain $x_j=t-u_j^{\frac{1}{n-\alpha}}$ and observe that $x_{k-1}<x_j<x_k$.  Thus, we obtain
\begin{eqnarray}\label{ResidDeriv}
R&=&\frac{1}{\Gamma(n+1-\alpha)}\sum_{k=1}^{m}\left(f^{(n)}(x_j)-\frac{f^{(n)}(x_k)+f^{(n)}(x_{k-1})}{2}\right)\left((t-x_{k-1})^{n-\alpha}-(t-x_k)^{n-\alpha}\right) \nonumber \\
&=&\frac{1}{2\Gamma(n+1-\alpha)}\sum_{k=1}^{m}\left((f^{(n)}(x_j)-f^{(n)}(x_k))+(f^{(n)}(x_j)-f^{(n)}(x_{k-1})\right)\left((t-x_{k-1})^{n-\alpha}-(t-x_k)^{n-\alpha}\right)\nonumber \\ 
&=&\frac{h_k}{2\Gamma(n+1-\alpha)}\sum_{k=1}^{m}\left(\frac{f^{(n)}(x_j)-f^{(n)}(x_k)}{h_k}+\frac{f^{(n)}(x_j)-f^{(n)}(x_{k-1})}{h_k}\right)\left((t-x_{k-1})^{n-\alpha}-(t-x_k)^{n-\alpha}\right)\nonumber \\
&=&\frac{h_k}{2\Gamma(n+1-\alpha)}\sum_{k=1}^{m}\left(-f^{(n+1)}(x_j)+O(h_k)+f^{(n+1)}(x_j)+O(h_k)\right)\left((t-x_{k-1})^{n-\alpha}-(t-x_k)^{n-\alpha}\right)\nonumber \\
&=&\frac{O(h_k^2)} {\Gamma(n+1-\alpha)}\sum_{k=1}^{m}\left((t-x_{k-1})^{n-\alpha}-(t-x_k)^{n-\alpha}\right)
\end{eqnarray}
Now, let us evaluate the difference 
\begin{eqnarray}
(t-x_{k-1})^{n-\alpha}-(t-x_k)^{n-\alpha} 
            &=&(t-x_{k}+h_k)^{n-\alpha}-(t-x_k)^{n-\alpha} \label{DeltaU} \\
            &=&  (p_k+h_k)^{n-\alpha}-p_k^{n-\alpha} < (p_k+h_{\max})^{n-\alpha}-p_k^{n-\alpha}. \nonumber
\end{eqnarray}
Here $p_k=t-x_k$.\\
The maximum of the last expression 
is achieved at $t-x_k=p_k=0$ since, in view of $n-\alpha-1 < 0$,
\begin{equation*}
	\frac{d}{dp}\left[(p+h_{\max})^{n-\alpha}-p^{n-\alpha}\right]\Big|_{p_k}=(n-\alpha)\left[(p_k+h_{\max})^{n-\alpha-1}-p_k^{n-\alpha-1}\right] <0. 
\end{equation*}
Hence, at $t-x_k=0$ expression \eqref{DeltaU} achieves its maximum which is less or equal $h_{\max}^{n-\alpha}$. Therefore 
\begin{equation*}
(t-x_{k-1})^{n-\alpha}-(t-x_k)^{n-\alpha} \le h_{\max}^{n-\alpha}, 
\end{equation*}
we obtain from \eqref{ResidDeriv}
\begin{equation*}
R \le O(h_{\max}^2)\cdot m\cdot h_{\max}^{n-\alpha} \le O(h_{\max}^{n-\alpha+1})C(t).
\end{equation*}
This proves Theorem \ref{Thm1}.
\end{proof}
%

%
\section{Substitution method for first-order linear equations}\label{SectionDiscr1}
In this section, we would like to analyze implementation of the numerical solution for the initial value problem of the linear fractional differential equation with $N$ fractional derivatives:
\begin{eqnarray}\label{LinFrDE}
 \sum_{l=1}^{N}q_l(t)D^{\alpha_l}y(t)+p(t)y(t)=f(t).
\end{eqnarray}
%
\begin{df}\label{DefOrderEqn}
Let $\alpha_{\max}=\max_{1\leq i\leq N}\{\alpha_i\}$. The order of fractional differential equation \eqref{LinFrDE} is $n\in \NN$, if $n-1<\alpha_{\max}<n$.
\end{df}
\noindent Let us consider equation \eqref{LinFrDE} of order one, $0<\alpha_l <1$.  In this case let the initial condition be 
\begin{equation}\label{LinFrDEIC} 
y(0)=y_0
\end{equation} 
and in view of \eqref{SubstDefFrD},
\begin{eqnarray}
D^{\alpha_l} y(t) = \ds\frac{1}{\Gamma(1-\alpha)}\int_{0}^{t}\frac{y'(x)dx}{(t-x)^{1-\alpha_l}}
\equiv\frac{1}{\Gamma(2-\alpha)}\int_{0}^{t^{n-\alpha}}y'(t-u^{\frac{1}{1-\alpha}})du. \nonumber
\end{eqnarray} 
 Partition of $x$ for all $D^{\alpha_l}y(x), l=1,...,N$ has to be the same to be able to apply finite differences: $x_0=0<x_1<...<x_m=t<...<x_M=hM$. Since equation \eqref{LinFrDE} can have derivatives of different order, it is important that $h=x_m-x_{m-1}=const, m=1,...,M$.
If we represent derivatives in form \eqref{FrDSubst}, then for each point $t=x_m=hm, m=1,...,M$ we get
\begin{eqnarray}\label{MainProblemDef}
\sum_{l=1}^{N}\frac{q_l(t)}{\Gamma(2-\alpha_l)}\sum_{k=1}^{m}\frac{y'(x_k)+y'(x_{k-1})}{2}\left((t-x_{k-1})^{1-\alpha_l}-(t-x_k)^{1-\alpha_l}\right)+p(t)y(t)=f(t). 
\end{eqnarray}
We use finite differences methods of the second order of approximation \cite{Eberly}:
\begin{eqnarray}\label{FDMs}
y'_k  &=&\frac{y_{k+1}-y_{k-1}}{2h}+O(h^2) \text { -- central difference,} \\
y'_k  &=&\frac{-3y_k+4y_{k+1}-y_{k+2}}{2h}+O(h^2) \text { -- forward difference,} \label{FDMf} \\
y'_k  &=&\frac{3y_k-4y_{k-1}+y_{k-2}}{2h}+O(h^2) \text { -- backward difference.} \label{FDMb}
\end{eqnarray}
We apply forward difference to the left end of the interval and backward differences to the right end.
Then,  equation \eqref{MainProblemDef} converts to the following equation for all points except $x=x_0=0$ ($y_0$ is given) as follows
\begin{eqnarray}\label{PointM_inEqn1}
&&\sum_{l=1}^{N}\frac{q_l(t)}{\Gamma(2-\alpha_l)} \Big[\frac{y'_1+y'_0}{2}\left((t-x_0)^{1-\alpha_l}-(t-x_1)^{1-\alpha_l}\right)+
\sum_{k=2}^{m-1}\frac{y'_k+y'_{k-1}}{2}\left((t-x_{k-1})^{1-\alpha_l}-(t-x_k)^{1-\alpha_l}\right) \nonumber \\
&+&\frac{y'_m+y'_{m-1}}{2}\left((t-x_{m-1})^{1-\alpha_l}-(t-x_m)^{1-\alpha_l}\right)\Big]
+p(t)y(t)=f(t).
\end{eqnarray}
Let us consider only the fractional derivative portion of equation \eqref{PointM_inEqn1}.
We put
\[
\phi_l(x_k)=(t-x_{k-1})^{1-\alpha_l}-(t-x_{k})^{1-\alpha_l}
\]
and approximate derivatives \eqref{FDMs}-\eqref{FDMb}. Then
\begin{eqnarray}\label{PointM_inEqn2}
&&\sum_{l=1}^{N}\frac{q_l(x_m)}{\Gamma(2-\alpha_l)} \Big[\frac{y_2-y_0 -3y_0+4y_1-y_2} {4h}\phi_l(x_1) 
+\sum_{k=2}^{m-1}\frac{y_{k+1}-y_{k-1}+y_{k}-y_{k-2}}{4h}\phi_l(x_k) \nonumber \\
&&\ +\frac{3y_m-4y_{m-1}+y_{m-2}+y_{m}-y_{m-2}}{4h}\phi_l(x_m)\Big] \nonumber \\
&&= \sum_{l=1}^{N}\frac{q_l(x_m)}{4h\Gamma(2-\alpha_l)} \Big[4(y_1-y_0) \phi_l(x_1) \nonumber \\
&&\ +\sum_{k=2}^{m-1}(-y_{k-2}-y_{k-1}+y_k+y_{k+1})\phi_l(x_k)
+4(y_m-y_{m-1})\phi_l(x_m)\Big].
\end{eqnarray}
Expression \eqref{PointM_inEqn2} can be presented as
\begin{equation}\label{PointM_inEqn4}
\sum_{l=1}^{N}\frac{q_l(x_m)}{4h\Gamma(2-\alpha_l)} \sum_{k=1}^{m}c^1_{lk}(x_0, x_1, ..., x_m)y_k,
\end{equation}
where $c^1_{lk}(x_0,x_1,...,x_m)$ are found as follows
\begin{subequations}\label{CalcCoeffsC}
\begin{align}
c^1_{l0} &= -4\phi_l(x_1)-\phi_l(x_2),\\
c^1_{l1} &= 4\phi_l(x_1)-\phi_l(x_2)-\phi_l(x_3),\\
c^1_{l2} &= \phi_l(x_1)+\phi_l(x_2)-\phi_l(x_3)-\phi_l(x_4),\\
...  \nonumber \\
c^1_{lk} &= \phi_l(x_{k-1})+\phi_l(x_k)-\phi_l(x_{k+1})-\phi_l(x_{k+2}),\\
... \nonumber \\
c^1_{l,m-1} &= \phi_l(x_{m-2})+\phi_l(x_{m-1})-4\phi_l(x_{m})),\\
c^1_{lm} &=\phi_l(x_{m-1})+4\phi_l(x_m).
\end{align}
\end{subequations}
Expression \eqref{PointM_inEqn4} can be rewritten as
\begin{equation*}\label{PointM_inEqn5}
\sum_{k=1}^{m}\left(\sum_{l=1}^{N}\frac{q_l(x_m)}{4h\Gamma(2-\alpha_l)} c^1_{lk}\right)\cdot y_k
\end{equation*}
and, assuming that
\begin{equation}\label{CalcCoeffsD}
 d_k=\sum_{l=1}^{N}\frac{q_l(x_m)}{4h\Gamma(2-\alpha_l)} c^1_{lk},
\end{equation}
we can represent equation \eqref{LinFrDE} as 
\begin{equation}\label{DiscrPres}
\sum_{k=0}^{m}d_k\cdot y_k + p_m y_m = f_m,
\end{equation}
where $p_m = p(x_m), f_m=f(x_m)$.
This, effectively, represents the discretization scheme of the substitution method to solve equation \eqref{LinFrDE}.  

\section{Second-order linear equations}\label{SectionDiscr2}
In the previous section we outlined the discretization of the first-order linear fractional differential equationsIn this setion we review second-order equations.\\
Let us consider linear equation \eqref{LinFrDE2} with $N$ derivatives
\begin{eqnarray}\label{LinFrDE2}
\sum_{l=1}^{s_1}q_l(t)D^{\alpha_l}y(t)+\sum_{s_1+1}^{N}q_l(t)D^{\alpha_l}y(t)+p(t)y(t)=f(t),\\
y(0)=y_0,\quad y'(0)=y'_0, \label{LinFrDE2IC}
\end{eqnarray}
where $0<\alpha_1<\alpha_2<...<\alpha_{s_1}<1<\alpha_{s_1+1}<...<\alpha_{N}<2$. 
Representation the derivatives in form \eqref{FrDSubst} yields
\begin{eqnarray}\label{MainProblemDef2}
&&\sum_{l=1}^{s_1}\frac{q_l(t)}{\Gamma(2-\alpha_l)}\sum_{k=1}^{m}\frac{y'(x_k)+y'(x_{k-1})}{2}\left((t-x_{k-1})^{1-\alpha_l}-(t-x_k)^{1-\alpha_l}\right) \phantom{WWWWWWWWWW}\nonumber \\
&&\ +\sum_{l=s_1+1}^{N}\frac{q_l(t)}{\Gamma(3-\alpha_l)}\sum_{k=1}^{m}\frac{y''(x_k)+y''(x_{k-1})}{2}\left((t-x_{k-1})^{2-\alpha_l}-(t-x_k)^{2-\alpha_l}\right) +p(t)y(t)=f(t).\phantom{W}
\end{eqnarray}
Here $x_m=t, x_m=hm, m=1,...,M$.  The first double sum corresponds to derivatives of order less than 1, the second double sum corresponds to the derivatives of order greater than 1 and less than 2.\\ 
We will use finite differences methods to represent the derivatives with the second order of approximation.  For the second derivatives we use the following formulas \cite{Eberly}:
\begin{eqnarray}\label{FDMs2}
y''_k  &=&\frac{y_{k-1}-2y_k+y_{k+1}}{h^2}+O(h^2) \text { -- central difference,} \\
y''_k  &=&\frac{2y_k-5y_{k+1}+4y_{k+2}-y_{k+3}}{h^2}+O(h^2) \text { -- forward difference,} \label{FDMf2} \\
y''_k  &=&\frac{2y_k-5y_{k-1}+4y_{k-2}-y_{k-3}}{2h}+O(h^2) \text { -- backward difference.} \label{FDMb2}
\end{eqnarray}
Assuming $\psi_l(x_k)=(t-x_{k-1})^{2-\alpha_l}-(t-x_{k})^{2-\alpha_l}$, equation \eqref{MainProblemDef2} converts to
\begin{eqnarray}\label{PointM_inEqn12}
&&\sum_{l=1}^{s_1}\frac{q_l(t)}{\Gamma(2-\alpha_l)} \Big[\frac{y'_1+y'_0}{2}\phi_l(x_1) 
+\sum_{k=2}^{m-1}\frac{y'_k+y'_{k-1}}{2}\phi_l(x_k) 
+\frac{y'_m+y'_{m-1}}{2}\phi_l(x_m)\Big] \label{PointM_inEqn12a} + p(t)y(t) \nonumber \\
&&\ +\hspace{-2mm} \sum_{l=s_1+1}^{N}\frac{q_l(t)}{\Gamma(3-\alpha_l)} \Big[\frac{y''_1+y''_0}{2}\psi_l(x_1)
+\sum_{k=2}^{m-1}\frac{y''_k+y''_{k-1}}{2}\psi_l(x_k)
+\frac{y''_m+y''_{m-1}}{2}\psi_l(x_m)\Big] =f(t).
\label{PointM_inEqn12c}
\end{eqnarray}
In the previous section in \eqref{PointM_inEqn2} we considered approximations for first-order derivatives. 
We approximate the second derivatives using \eqref{FDMs2}-\eqref{FDMb2} as follows:
\begin{eqnarray}\label{PointM_inEqn22}
&&\sum_{l=s_1+1}^{N}\frac{q_l(x_m)}{\Gamma(3-\alpha_l)} \Big[\frac{y_2-2y_1+y_0+2y_0-5y_1+4y_2-y_3} {2h^2}\psi_l(x_1) \nonumber \\
&&\ +\sum_{k=2}^{m-1}\frac{y_{k+1}-2y_k+y_{k-1}+y_{k}-2y_{k-1}+y_{k-2}}{2h^2}\psi_l(x_k) \nonumber \\
&&\ +\frac{2y_m-5y_{m-1}+4y_{m-2}-y_{m-3}+y_{m}-2y_{m-1}+y_{m-2}}{2h^2}\psi_l(x_m) \Big] \nonumber \\
&=& \sum_{l=s_1+1}^{N}\frac{q_l(x_m)}{2h^2\Gamma(3-\alpha_l)} \Big[(3y_0-7y_1+5y_2-y_3) \psi_l(x_1) \nonumber \\
&&\ +\sum_{k=2}^{m-1}(y_{k-2}-y_{k-1}-y_k+y_{k+1})\psi_l(x_k) 
+(3y_m-7y_{m-1}+5y_{m-2}-y_{m-3})\psi_l(x_m) \Big].
\end{eqnarray}
Expression \eqref{PointM_inEqn22} can be presented as
\begin{equation}\label{PointM_inEqn23}
\sum_{l=s_1+1}^{N}\frac{q_l(x_m)}{2h^2\Gamma(3-\alpha_l)} \sum_{k=1}^{m}c^2_{lk}(x_0, x_1, ..., x_m)y_k,
\end{equation}
where $c^2_{lk}(x_0,x_1,...,x_m)\equiv c^2_{lk}$ can be found as follows
\begin{subequations}\label{CalcCoeffsC2}
\begin{align}
c^2_{l0} &= 3\psi_l(x_1)+\psi_l(x_2),\\
c^2_{l1} &= -7\psi_l(x_1)-\psi_l(x_2)+\psi_l(x_3),\\
c^2_{l2} &= 5\psi_l(x_1)-\psi_l(x_2)-\psi_l(x_3)+\psi_l(x_4),\\
c^2_{l3} &= -\psi_l(x_1)+\psi_l(x_2)-\psi_l(x_3)-\psi_l(x_4)+\psi_l(x_5),\\
...  \nonumber \\
c^2_{lk} &= \psi_l(x_{k-1})-\psi_l(x_k)-\psi_l(x_{k+1})+\psi_l(x_{k+2}),\\
... \nonumber \\
c^2_{l,m-3} &= \psi_l(x_{m-4})-\psi_l(x_{m-3})-\psi_l(x_{m-2}+\psi_l(x_{m-1}-\psi_l(x_m),\\
c^2_{l,m-2} &= \psi_l(x_{m-3})-\psi_l(x_{m-2})-\psi_l(x_{m-1}+5\psi_l(x_m),\\
c^2_{l,m-1} &= \psi_l(x_{m-2})-\psi_l(x_{m-1})-7\psi_l(x_{m}),\\
c^2_{lm} &=3\psi_l(x_{m-1})+\psi_l(x_m).
\end{align}
\end{subequations}
Now, we can re-write expression \eqref{PointM_inEqn23} as
\begin{equation*}\label{PointM_inEqn43}
\sum_{k=1}^{m}\left(\sum_{l=s_1+1}^{N}\frac{q_l(x_m)}{2h^2\Gamma(3-\alpha_l)} c^2_{lk}\right) y_k
\end{equation*}
and, therefore, equation \eqref{PointM_inEqn12} becomes
\begin{equation*}
\sum_{k=1}^{m}\left[\frac{1}{2h^2}\left(\sum_{l=1}^{s_1}\frac{q_l(x_m)}{\Gamma(2-\alpha_l)} c^1_{lk} +
\sum_{l=s_1+1}^{N}\frac{q_l(x_m)}{\Gamma(3-\alpha_l)} c^2_{lk}\right)\right] y_k +p(x_m)u(x_m)=f(x_m).
\end{equation*}
Assuming that 
\begin{equation}\label{CalcCoeffsD2}
d_k=\frac{1}{2h^2}\left(\sum_{l=1}^{s_1}\frac{q_l(x_m)}{\Gamma(2-\alpha_l)} c^1_{lk} +
\sum_{l=s_1+1}^{N}\frac{q_l(x_m)}{\Gamma(3-\alpha_l)} c^2_{lk}\right),
\end{equation}
Equation \eqref{LinFrDE2} can be represented as 
\begin{equation}\label{DiscrPres2}
\sum_{k=0}^{m}d_k\cdot y_k + p_m y_m = f_m,
\end{equation}
where $p_m = p(x_m), f_m=f(x_m)$.  This equation is the same as equation \eqref{DiscrPres} except for the values of $d_k$.
This represents the discretization scheme of the substitution method to solve the second-order 
equation \eqref{LinFrDE2}. 

\section{Higher-order equations}\label{SectionDiscrN}
In the previous two sections we developed the discretization schemes for linear fractional differential equations for $0<\alpha<2$.  If equation \eqref{LinFrDE2} contains derivatives of order higher than two, then the derivatives of higher order will need to be represented through finite differences similar to 
\eqref{MainProblemDef2}-\eqref{FDMb2}.  
We consider equation of order $r$
\begin{eqnarray}\label{LinFrDEk}
&&\sum_{n=1}^{r}\left(\sum_{l=s_{n-1}+1}^{s_n}q_l(t)D^{\alpha_l}y(t)\right)+p(t)y(t)=f(t),\\
&&y(0)=y_0, y'(0)=y'_0, y''(0)=y''_0,..., y^{(r-1)}(0)=y^{(r-1)}_0. \label{LinFrDEkIC}
\end{eqnarray}
Here $0<\alpha_l<r, l=1,...,N$ and the integer values of $s_n$ show how many fractional derivatives fall into interval $(n-1,n)$: $0<\alpha_1<...<\alpha_{s_1}<1<\alpha_{s_1+1}<...<\alpha_{s_2}<2<...<\alpha_{s_{r-1}}<r-1<...<\alpha_{s_r}=\alpha_N<r$.  $N$ is the total number of fractional derivatives and therefore   
$s_0=0<s_1<s_2<...<s_r=N$.\\ 

Using the substitution approximation, we can write equation \eqref{LinFrDEk} as
\begin{eqnarray}\label{MainProblemDefk}
&&\sum_{n=1}^{r}\left(\sum_{l=s_{n-1}+1}^{s_n}\frac{q_l(t)}{\Gamma(n+1-\alpha_l)}\sum_{k=1}^{m}\frac{y^{(n)}(x_k)+y^{(n)}(x_{k-1})}{2}\left((t-x_{k-1})^{n-\alpha_l}-(t-x_k)^{n-\alpha_l}\right)\right)\nonumber \\
&+&p(t)y(t)=f(t).
\end{eqnarray}
Here $x_m=t, x_m=hm, m=1,...,M$. 

The approach and basic logic of the implementation remain the same.  The second order of approximation  of the $n$-order derivative is  
\begin{eqnarray}\label{FDMsn}
y^{(n)}_k  &=&\frac{1}{B_n h^n}\sum_{l=k-n_2}^{k+n_2} a_l y_l +O(h^2) \text { -- central difference,} \\
y^{(n)}_k  &=&\frac{1}{B_n h^n}\sum_{l=k}^{k+n+1} a_l y_l+O(h^2) \text { -- forward difference,} \label{FDMfn} \\
y^{(n)}_k  &=&\frac{1}{B_n h^n}\sum_{l=k-n-1}^{k} a_l y_l +O(h^2) \text { -- backward difference.} \label{FDMbn}
\end{eqnarray}
where $B_n = \begin{cases}
1, & \text{if}\ n \text { is even} \\
2, & \text{if}\ n \text { is odd}\end{cases}$; $ n_2 = \left\lceil \frac{n}{2}\right\rceil$ $\left(\frac{n}{2} \text{ rounded up to the nearest integer}\right)$ and coefficients $a_l$ can be calculated as a solution to the linear equation system.

In order to figure out the coefficients $a_l$ in \eqref{FDMsn}, we consider the   
\textbf{central difference} coefficients system for the $n^{th}$ derivative:

$\begin{pmatrix}
1 & 1 & 1 & ... & 1 \\
-n_2 & -n_2+1 & ... & n_2-1 & n_2 \\
\frac{n_2^2}{2} & \frac{(n_2-1)^2}{2}& ... & \frac{(n_2-1)^2}{2} & \frac{n_2^2}{2} \\
... & ... & ... & ... & ...\\
\frac{(-n_2)^l}{l!} & \frac{(-n_2+1)^{l}}{l!}  & ... &\frac{(n_2-1)^l}{l!} & \frac{n_2^l}{l!}\\ 
... & ... & ... & ... & ...\\
\frac{(-n_2)^n}{n!} & \frac{(-n_2+1)^{n}}{n!}  & ... &\frac{(n_2-1)^n}{n!} & \frac{n_2^n}{n!}\\ 
\frac{(-n_2)^{n+1}}{(n+1)!} & \frac{(-n_2+1)^{n+1}}{(n+1)!}  & ... & \frac{(n_2-1)^{n+1}}{(n+1)!}& \frac{n_2^{n+1}}{(n+1)!} 
\end{pmatrix}$
$\begin{pmatrix}
a_{k-n_2} \\ a_{k-n_2+1}\\ a_{k-n_2+2}\\ ... \\ a_l \\ ... \\a_{k+n_2-1} \\ a_{k+n_2}
\end{pmatrix}$
=
$\begin{pmatrix}
0 \\ 0 \\ 0 \\ ... \\ 0 \\ ...\\ 1 \\ 0
\end{pmatrix}$.

It is based on summation \eqref{FDMsn}, where $y_l=y(x_k+(l-k)h)$ is represented through Taylor series as shown in \cite{Eberly}. The coefficients $a_l$ can be verified from the algorithm and the tables developed in \cite{Fornberg}.

If we look at the $n$-th derivative as $(y_k^{(n-2)})''$ then the central difference of the second order can be recursively expressed as 
\begin{equation}\label{nthderivassecond}
y_k^{(n)} = \frac{y_{k-1}^{(n-2)}-2y_k^{(n-2)}+y_{k+1}^{(n-2)}}{h^2}. 
\end{equation}
As we see here, if the second derivative is applied to the two-less derivative $y^{n-2}$, it expands that derivative by two elements. However, if we apply the first derivative twice, then the expression is expanded by 4 elements.
Certainly, chosen precision ($O(h^2)$ in our case), we try to minimize the number of elements of $\{y_k\}$.
%
\begin{thm}\label{Thm2}
	Coefficients in the denominator of central differences of the second order are the same as the coefficients of  polynomial $(a-b)^n$ (binomial coefficients) for even derivatives ($n$ is even) and coincide with the coefficients of the polynomial \\ $(a-b)^{n-1}(a^2-b^2)$ for odd derivatives ($n$ is odd). 
\end{thm}
\begin{proof}
	Let's address the first statement of the theorem -- the coefficients for even derivatives represented through central differences of the second order.  This can be proved by induction.  We know that second derivative is represented as $y''_k  =\ds\frac{y_{k-1}-2y_k+y_{k+1}}{h^2}$.  We also know that $y_k^{(4)}=\left(y_k''\right)''$, therefore we get 
	\begin{eqnarray}\label{4thorderCoeffs}
	y_k^{(4)}&=&\frac{1}{h^2}(y''_{k-1}-2y''_k+y''_{k+1})\nonumber \\
	&=&\frac{1}{h^4}(y_{k-2}-2y_{k-1}+y_k-2y_{k-1}+4y_k-2y_{k+1}+y_{k+2}-2y_{k+1}+y_k) \nonumber \\
	&=&\frac{1}{h^4}(y_{k-2}-4y_{k-1}+6y_k-4y_{k+1}+y_{k+2}).
	\end{eqnarray}
	Coefficients in parenthesis of \eqref{4thorderCoeffs} are exactly the same as coefficients of the polynomial $(x-y)^4$. \\
	Now, let's assume that it's true for the $n_{th}$ even derivative:
	\begin{equation}\label{nthEvenDeriv}
	y_l^{(n)} = \sum_{k=0}^{n}\binom{n}{k}(-1)^k y_{l-\frac{n}{2}+k}
	\end{equation}
	Let us verify that the formula is correct for $n+2$.
	\begin{eqnarray}\label{np2EvenDeriv}
	y_l^{(n+2)} &=&\left(y_l^{(n)}\right)''= \frac{1}{h^2}\left(y_{l+1}^{(n)}-2y_l^{(n)}+y_{l-1}^{(n)}\right)\nonumber \\
	&=&\frac{1}{h^4}\left(\sum_{k=0}^{n}\binom{n}{k}(-1)^k y_{l+1-\frac{n}{2}+k}-2\sum_{k=0}^{n}\binom{n}{k}(-1)^k y_{l-\frac{n}{2}+k}+\sum_{k=0}^{n}\binom{n}{k}(-1)^k y_{l-1-\frac{n}{2}+k}\right)\nonumber \\
	&=&\frac{1}{h^4} \sum_{k=0}^{n}\binom{n}{k}(-1)^k \left(y_{l+1-\frac{n}{2}+k}-2y_{l-\frac{n}{2}+k} +y_{l-1-\frac{n}{2}+k}\right).
	\end{eqnarray}
	Coefficients of the first and the last terms in expression \eqref{np2EvenDeriv} are equal to one.  For the second and one before the last coefficients the are
	\begin{eqnarray}\label{binom2}
	-\binom{n}{1}-2\binom{n}{0}&=&-n-2=-(n+2)	\text { -- second term},\\  
	-2\binom{n}{n}-\binom{n}{n-1}&=&-2-n=-(n+2) \text {-- one before last term}.
	\end{eqnarray}
	Finally, let's check the rest of the coefficients. The expression for the $k_{th}$ element can be simplified:
	\begin{equation}\label{binomproof}
	(-1)^k \left[\binom{n}{k}+2\binom{n}{k-1}+\binom{n}{k-2}\right]=(-1)^k {n+2 \choose k}.
	\end{equation}
	Formulas \eqref{binom2}-\eqref{binomproof} prove that for even derivatives the statement of the theorem is correct. \\
	The odd derivative can be calculated either based on the statement \eqref{nthderivassecond}, applied to the odd derivative (the proof is almost identical to the proof for even derivatives) or by considering the first derivative of the previous even derivative:
	\begin{equation}\label{nthderivasfirst}
	y_k^{(n)}=\frac{y_{k+1}^{(n-1)}-y_{k-1}^{(n-1)}}{2h}, \quad n \mbox{\ is odd}.
	\end{equation}
	Since $n-1$ is even, then each derivative in \eqref{nthderivasfirst} can be represented as a sum of elements of $\{y_k\}$ with binomial coefficients like in $(a-b)^{n-1}$.  If we subtract $(n-1)^{st}$ odd derivative from another of two elements before -- $(n+1)^{st}$, it is equivalent to shifting of all elements of the latter to the right and subtracting them.  This action can be represented as multiplication of $(a-b)^{n-1}$ by $(a^2-b^2)$ and extracting the coefficient in front of all different terms.  The first two terms are the terms of the first polynomial, the last two -- of the second one and all the ones in the middle are subtracted from one another. \\
	This proves Theorem \ref{Thm2}.
\end{proof}

\noindent Below is an example of representation of the shift in Theorem \ref{Thm2}.  It is for the polynomial of the 4-th degree\\
$\begin{matrix}
1a^6 & -4a^5b & 6a^4b^2  & -4a^3b^3 & 1a^2b^4  &   & \\
&        &-1b^2a^4  & 4b^3a^3  & -6b^4a^2 & 4b^5a & -1b^6\\
----  &----   &----&----&----&----&----\\
1a^6 & -4a^5b & 5a^4b^2 &  0b^3a^3  & -5a^2b^4 & 4b^5a & -1b^6.	 
\end{matrix} $\\
It presents the calculation of the coefficients for the 5th derivative through central difference. 

\textbf{Remark}: if we look at application of the second derivative, we can see that from the point of view of the coefficients, it represents their summation shifted to the left with the set shifted to the right minus twice the actual set therefore it is equivalent to $(a-b)^n$ multiplied by $a^2+b^2-2ab$ which gives us $(a-b)^{n+2}$, which is the statement of the first part of the theorem.  

In order to identify coefficients $a_l$ in \eqref{FDMfn}, \eqref{FDMbn} we must consider
\textbf{Forward difference} coefficients system for the $n^{th}$ derivative

$\begin{pmatrix}
1 & 1 & 1 & ... & 1 \\
0 & 1 & 2 & ... & n+1 \\
... & ... & ... & ... & ...\\
0 & \frac{1}{l!} & \frac{2^{l}}{l!}  & ... & \frac{(n+1)^l}{l!}\\ 
... & ... & ... & ... & ...\\
0 & \frac{1}{n!} & \frac{2^{n}}{n!}  & ... & \frac{(n+1)^n}{n!}\\ 
0 & \frac{1}{(n+1)!} & \frac{2^{n+1}}{(n+1)!}  & ... & \frac{(n+1)^{n+1}}{(n+1)!} 
\end{pmatrix}$
$\begin{pmatrix}
a_k \\ a_{k+1}\\
... \\ a_l \\ ... \\ a_{k+n} \\ a_{k+n+1}
\end{pmatrix}$
=
$\begin{pmatrix}
0 \\ 0 \\
... \\ 0 \\ ... \\ 1 \\ 0
\end{pmatrix}$ 

and \textbf{Backward difference} coefficients system for the $n^{th}$ derivative:

$\begin{pmatrix}
1 & 1 & 1 & ... & 1 \\
0 & 1 & 2 & ... & n+1 \\
... & ... & ... & ... & ...\\
0 & \frac{1}{l!} & \frac{2^{l}}{l!}  & ... & \frac{(n+1)^l}{l!}\\ 
... & ... & ... & ... & ...\\
0 & \frac{1}{n!} & \frac{2^{n}}{n!}  & ... & \frac{(n+1)^n}{n!}\\ 
0 & \frac{1}{(n+1)!} & \frac{2^{n+1}}{(n+1)!}  & ... & \frac{(n+1)^{n+1}}{(n+1)!} 
\end{pmatrix}$
$\begin{pmatrix}
a_k \\ a_{k-1}\\ 
... \\ a_l \\ ... \\ a_{k-n} \\ a_{k-n-1}
\end{pmatrix}$
=
$\begin{pmatrix}
0 \\ 0 \\
... \\0 \\...\\ (-1)^n \\ 0
\end{pmatrix}$.

As we can see there is a set of terms that is added to in equation \eqref{MainProblemDefk} for each range $n-1<\alpha_l<n$:
\begin{equation}\label{DerivRange}
\sum_{l=s_{n-1}+1}^{s_n}\frac{q_l(t)}{\Gamma(n+1-\alpha_l)}\sum_{k=1}^{m}\frac{y^{(n)}(x_k)+y^{(n)}(x_{k-1})}{2}\left((t-x_{k-1})^{n-\alpha_l}-(t-x_k)^{n-\alpha_l}\right),
\end{equation}
Taking into consideration the fact that central finite differences of the second order are represented with $n_2$ elements prior and after $y_k$ for $y_k^{(n)}$ as in \eqref{FDMsn}, there will have to be $n_2$ number of forward and backward differences used in expression \eqref{PointM_inEqn22}.  In order to find $n^{th}$ derivatives for the first $n_2$ points on the interval we need to use forward difference formulas, for the last $n_2$ points -- backward differences.  The derivatives for the rest of the points on the grid are sought through central differences. 

We can apply formulas \eqref{CalcCoeffsC} for the derivatives of the first order, formulas \eqref{CalcCoeffsC2} for the second order and derive similar expressions for coefficients of each order $r$ \eqref{DerivRange} as was done in \eqref{CalcCoeffsC} and \eqref{CalcCoeffsC2} for \eqref{PointM_inEqn4} and \eqref{PointM_inEqn23}.  Then, expression \eqref{DerivRange} becomes
\begin{equation*}\label{PointM_inEqn4k}
\sum_{k=1}^{m}\left(\sum_{l=s_{n-1}+1}^{s_n}\frac{q_l(x_m)}{2h^2\Gamma(n+1-\alpha_l)} c^n_{lk}\right) y_k
\end{equation*}
and therefore equation \eqref{MainProblemDefk} takes the following form
\begin{equation*}
\sum_{k=1}^{m}\left[\frac{1}{2h^2}\sum_{n=1}^{r}\left(\sum_{l=s_{n-1}+1}^{s_n}\frac{q_l(x_m)}{\Gamma(n+1-\alpha_l)} c^n_{lk}\right)\right] y_k +p(x_m)u(x_m)=f(x_m).
\end{equation*}
Assuming that 
\begin{equation}\label{CalcCoeffsDk}
d_k=\frac{1}{2h^2}\sum_{n=1}^{r}\left(\sum_{l=s_{n-1}+1}^{s_n}\frac{q_l(x_m)}{\Gamma(n+1-\alpha_l)} c^n_{lk}\right),
\end{equation}
Equation \eqref{LinFrDEk} of order $r$ can be represented as 
\begin{equation}\label{DiscrPresk}
\sum_{k=0}^{m}d_k\cdot y_k + p_m y_m = f_m,
\end{equation}
where $p_m = p(x_m), f_m=f(x_m)$.  This equation is the same as equations \eqref{DiscrPres}, \eqref{DiscrPres2} except for the values of $d_k$.

\section{Sufficiency for linear system to be well-conditioned}\label{SectionSuffCond}
It is important to emphasize that equations \eqref{DiscrPres}, \eqref{DiscrPres2} and \eqref{DiscrPresk} are identical.  The only differences are the definition of the coefficients $d_i$ as reflected in expressions \eqref{CalcCoeffsD}, \eqref{CalcCoeffsD2} and \eqref{CalcCoeffsDk} and the number of initial conditions \eqref{LinFrDEIC}, \eqref{LinFrDE2IC}, \eqref{LinFrDEkIC}.  For simplicity we analyze the first order case corresponding to equation \eqref{DiscrPres}:
\begin{eqnarray}\label{MainProblemThrm}
\sum_{k=0}^{m-1}d_k\cdot y_k + (d_m+p_m) y_m  &=& f_m, 1\le m \le M,\\
y_0&=& \mu. \label{MainProblemThrm2}
\end{eqnarray} 
\begin{df}{}
We say that the difference initial value problem \eqref{MainProblemThrm}, \eqref{MainProblemThrm2} with bounded coefficients $|p_m|, |d_k| < P,  k=0,1,...m,$ is well-conditioned if for all large enough $M$ it has only one solution, $\{y_m\}$, for arbitrary right hand sides $\mu$ and $\{f_m\}$, and if the numerical solution $\{y_m\}$ satisfies the bound 
\begin{equation*}
|y_m| \le K \max\{|\mu|, \max\limits_{m}|f_m|\},
\end{equation*}   
where $K>0$ is a constant which does not depend on $M$.
\end{df} 
%
\begin{thm}\label{Thrm3}
If coefficients $d_k$ and $p_m$ in \eqref{MainProblemThrm}  satisfy for some $\delta>0$ the condition
\begin{equation}\label{Thrm2Statement}
|d_m+p_m| \ge \sum_{k=0}^{m-1}|d_k| + \delta
\end{equation}
for all $m=1,...,M$, then problem \eqref{MainProblemThrm}, \eqref{MainProblemThrm2} is well conditioned and the solution $\{y_m\}$ satisfies the following upper estimate
\begin{equation}\label{AddlResThrm2}
|y_m|\le \max \{|\mu|, \frac{1}{\delta}\max_{k=0,1,...,m}|f_k| \}.
\end{equation}
\end{thm}
%
\begin{proof}{}
First let's assume that the problem \eqref{MainProblemThrm}, \eqref{MainProblemThrm2} for given $\mu, \{f_m\}$ has a solution $\{y_m\}$ and verify that this solution satisfies \eqref{AddlResThrm2}.  Let $y_L=\max_m \{|y_m|\}$.  For $y_L=y_0$, \eqref{AddlResThrm2} is obvious since $y_0=\mu$.  Now, let's consider cases $0<L\le M, |y_L| \ge |y_m|, 0\le m\le M$.  Then, taking into account \eqref{MainProblemThrm}, we can state
%
\bq	|d_L+p_L|\cdot |y_L| = |-\sum_{k=0}^{L-1}d_k y_k + f_L| 
	\le \sum_{k=0}^{L-1}|d_k|\cdot |y_k| + |f_L| \le |y_L|\sum_{k=0}^{L-1}|d_k| + |f_L|.
\eq
Therefore
\[
|y_m|\le |y_L| \le \frac{|f_L|}{|d_L+p_L| - \sum_{k=0}^{L-1}|d_k|} \le \frac{|f_L|}{\delta},
\]
and we obtain that \eqref{AddlResThrm2} is also satisfied.  \\
Now we need to show that problem \eqref{MainProblemThrm}, \eqref{MainProblemThrm2} has only one solution $\{y_m\}$ for any given right-hand sides $\mu, \{f_m\}$.  The given problem can be viewed as a system of $M+1$ linear equations with $M+1$ unknowns $\{y_m\}_{m=0}^{m=M}$. The determinant of the system of equations is not equal to zero if the corresponding homogeneous system has only a trivial solution ($y_m=0, m=0,1,...,M$).  In our problem, homogeneous system will have $\mu = 0, f_m=0$ for all $m,0\le m \le M$.  From inequality \eqref{AddlResThrm2} we 
obtain $|y_m| \equiv 0$ and there is only trivial solution for the homogeneous system of equations. 
This proves Theorem \ref{Thrm3}. 
\end{proof}

\noindent\textbf{Remark}. We can add that there is another sufficient condition for well-conditioning of the problem:
\begin{eqnarray}\label{AddlWPCond}
	\frac{|d_m+p_m|-\sum_{k=0}^{m-1}|d_k|}{|d_m+p_m|+\sum_{k=0}^{m-1}|d_k|} \ge A > 0, \qquad \qquad 
	\max_m \{|d_m+p_m|, \sum_{k=0}^{m-1}|d_k|\} \ge B >0,
\end{eqnarray}
where $A$ and $B$ are positive constants independent of $M$ and $m$.  From condition \eqref{AddlWPCond} we get condition \eqref{Thrm2Statement} with the constant $A$:
\begin{equation*}
\delta = A \left(|d_m+p_m|+\sum_{k=0}^{m-1}|d_k|\right) \ge A\cdot B > 0,
\end{equation*} 
then \eqref{AddlResThrm2} takes the form
\begin{equation*}
|y_m| \le \max \{|\mu|, \frac{1}{A\cdot B}\max_{k=0,1,...,m}|f_k| \}.
\end{equation*}

In the case of the second order equation with initial conditions \eqref{LinFrDE2IC}, the second equation in the problem \eqref{MainProblemThrm}, \eqref{MainProblemThrm2} becomes: $y_0= \mu; \quad y_1= \nu$, because $y'(0)=y'_0$ can be represented through finite differences as $y'_0=\ds\frac{y_1-y_0}{h}+O(h)$, and therefore $y_1=hy'_0+y_0=h y'_0+\mu \equiv \nu$.  In this case, we need to replace $\mu$ in \eqref{AddlResThrm2} with $\max(\mu,\nu)$ and the rest of the sufficient condition stays the same (point $y_1$ is satisfied and does not need to be checked).  The theorem does not change.  

If equation is of higher order $n$ \eqref{LinFrDEk}, then $n$ initial conditions need to be provided.  Similar logic can be used to identify the first $n$ elements of the solution and the $\max_{k=0,1,...,n}|y_k|$ should be used in the theorem as $\mu$. 

On the surface it looks like these conditions for the well-conditioned system are very difficult to check, but in reality, they are easily implemented as a program for any equation of type \eqref{LinFrDE}, \eqref{LinFrDE2} or \eqref{LinFrDEk} through the use of the formulas \eqref{CalcCoeffsD}, \eqref{CalcCoeffsD2}, \eqref{CalcCoeffsDk}, and taking into account Theorem \ref{Thrm3}.  This theorem represents the sufficient condition for the system to be well-conditioned.  

\section{Examples of a problem solved by numerical implementation of the substitution method}\label{SectionExamples}
\begin{figure}[H]
	\begin{center}
		\includegraphics[width=7.5cm]{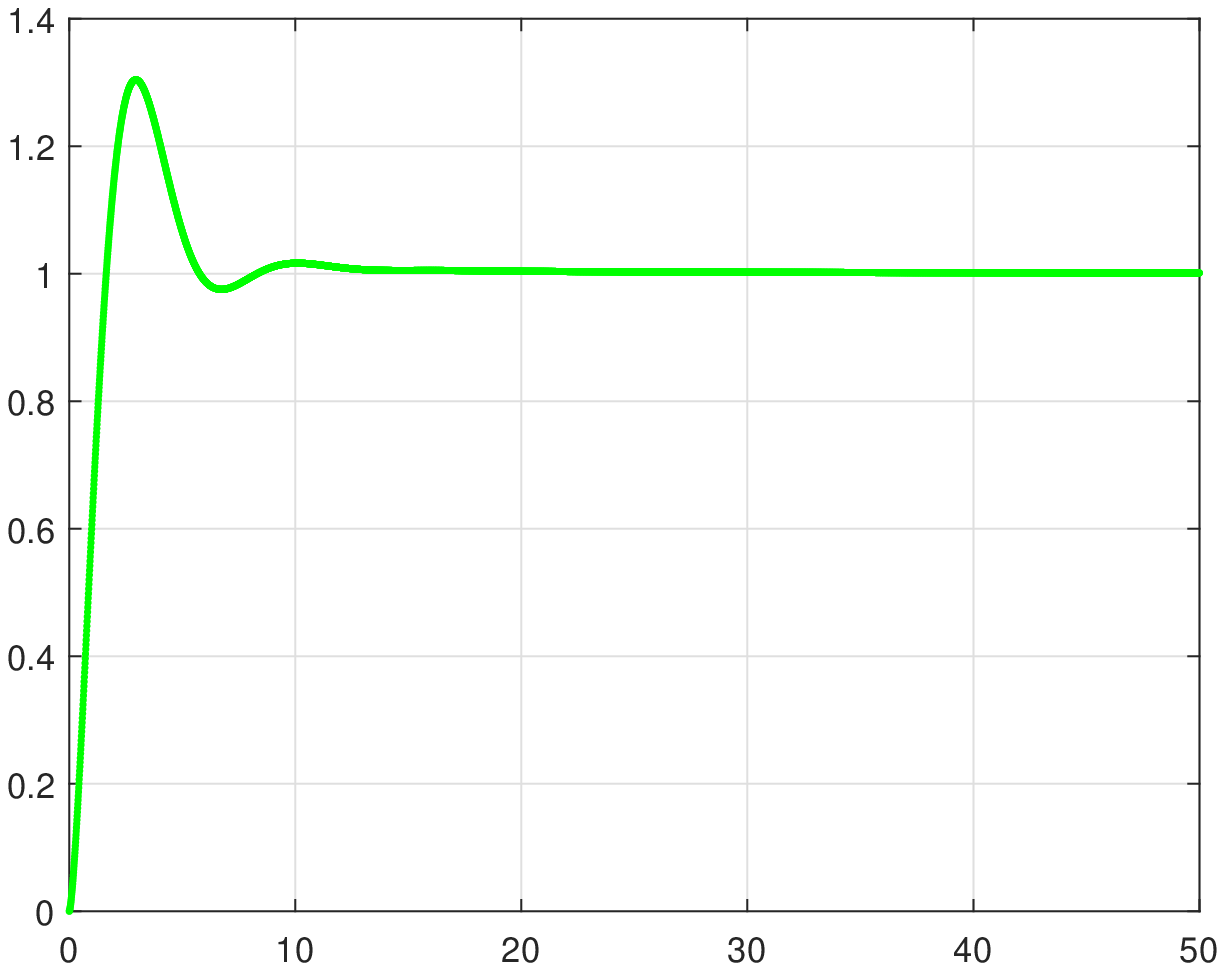}
		\includegraphics[width=7.5cm]{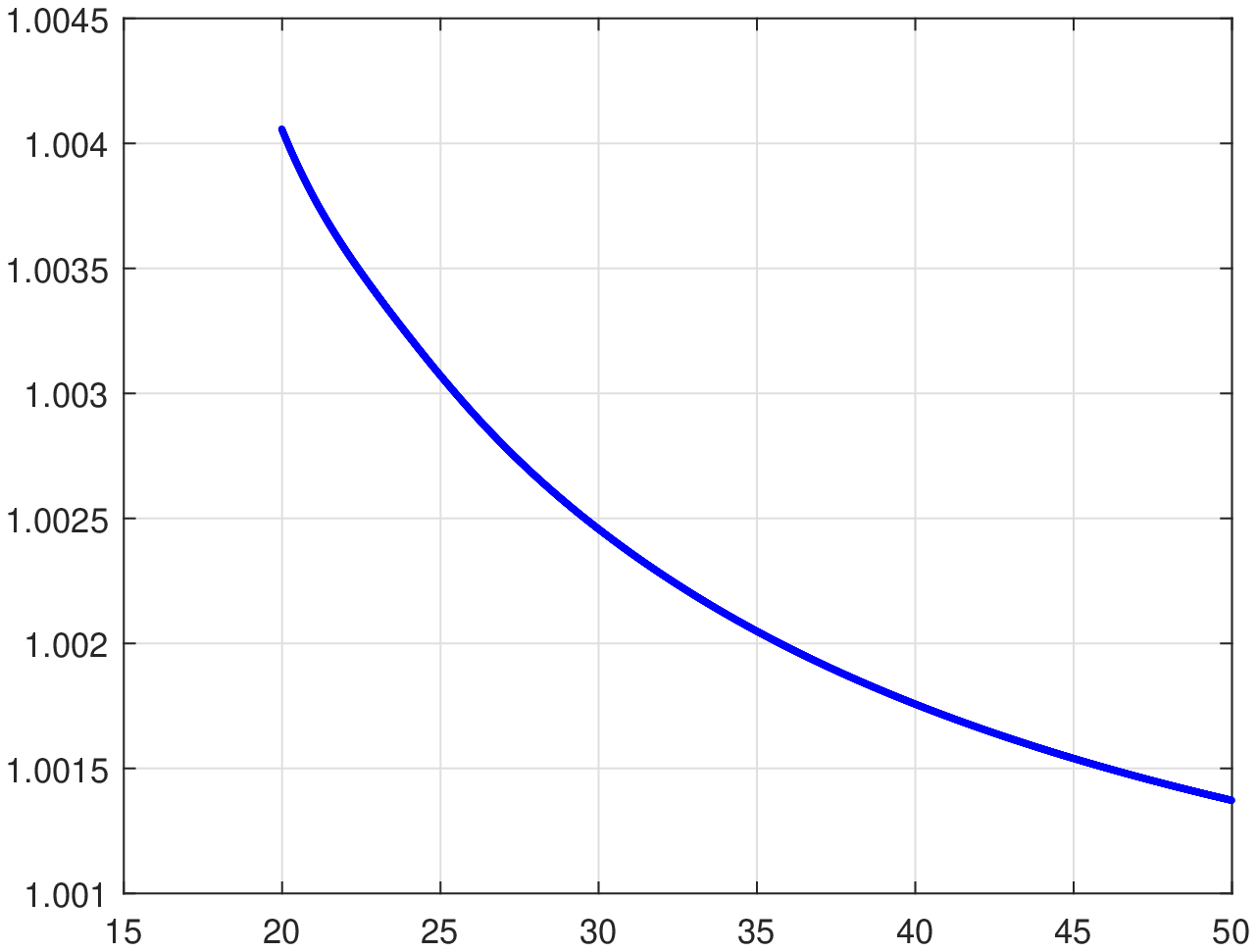}
		\caption{$f(x)=1$. Left: full solution, Right: zoom -- matches \cite{Podlubny}.}\label{ExampleFig1}
	\end{center}
\end{figure}

\begin{figure}[H]
	\begin{center}
		\includegraphics[width=7.5cm]{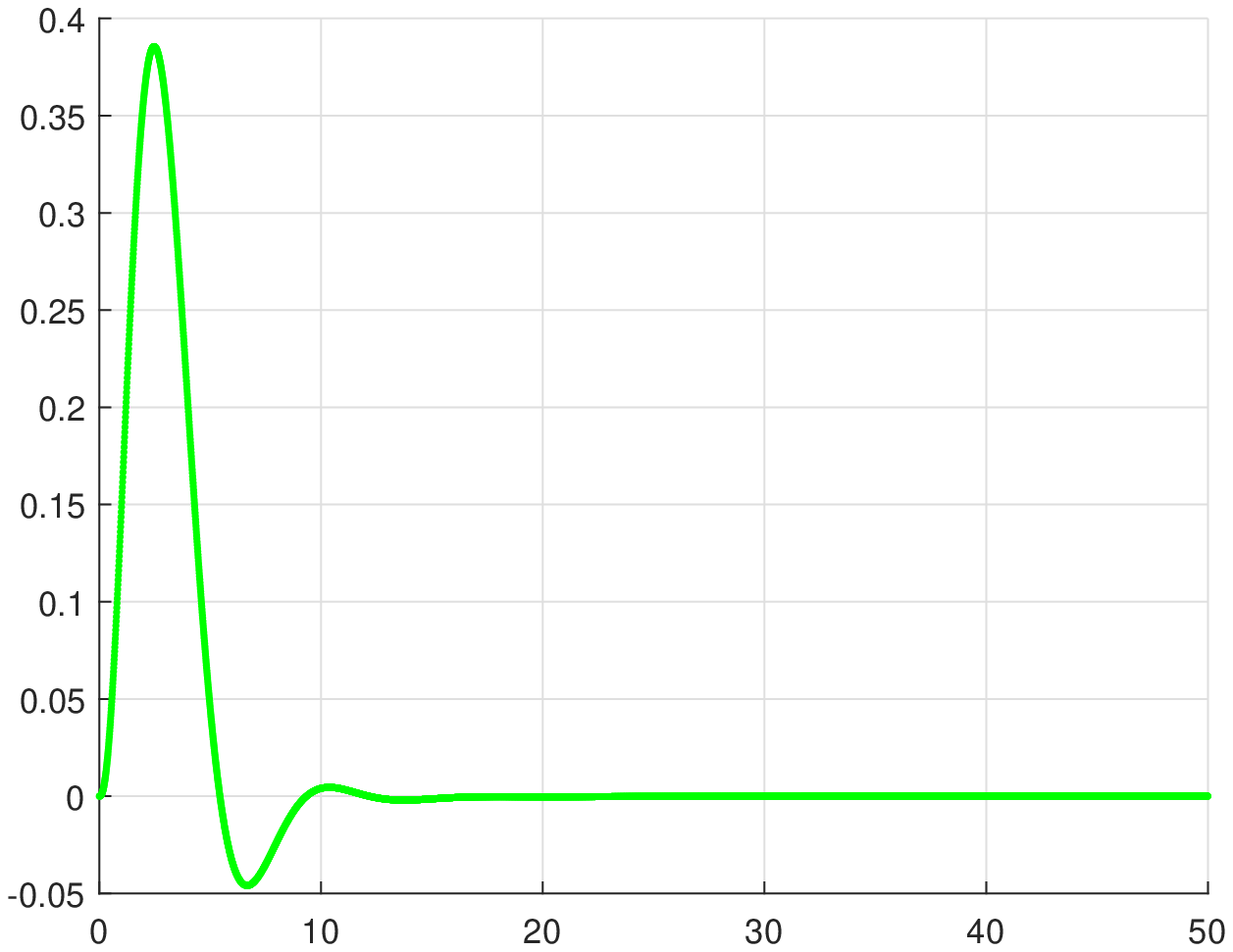}
		\includegraphics[width=7.5cm]{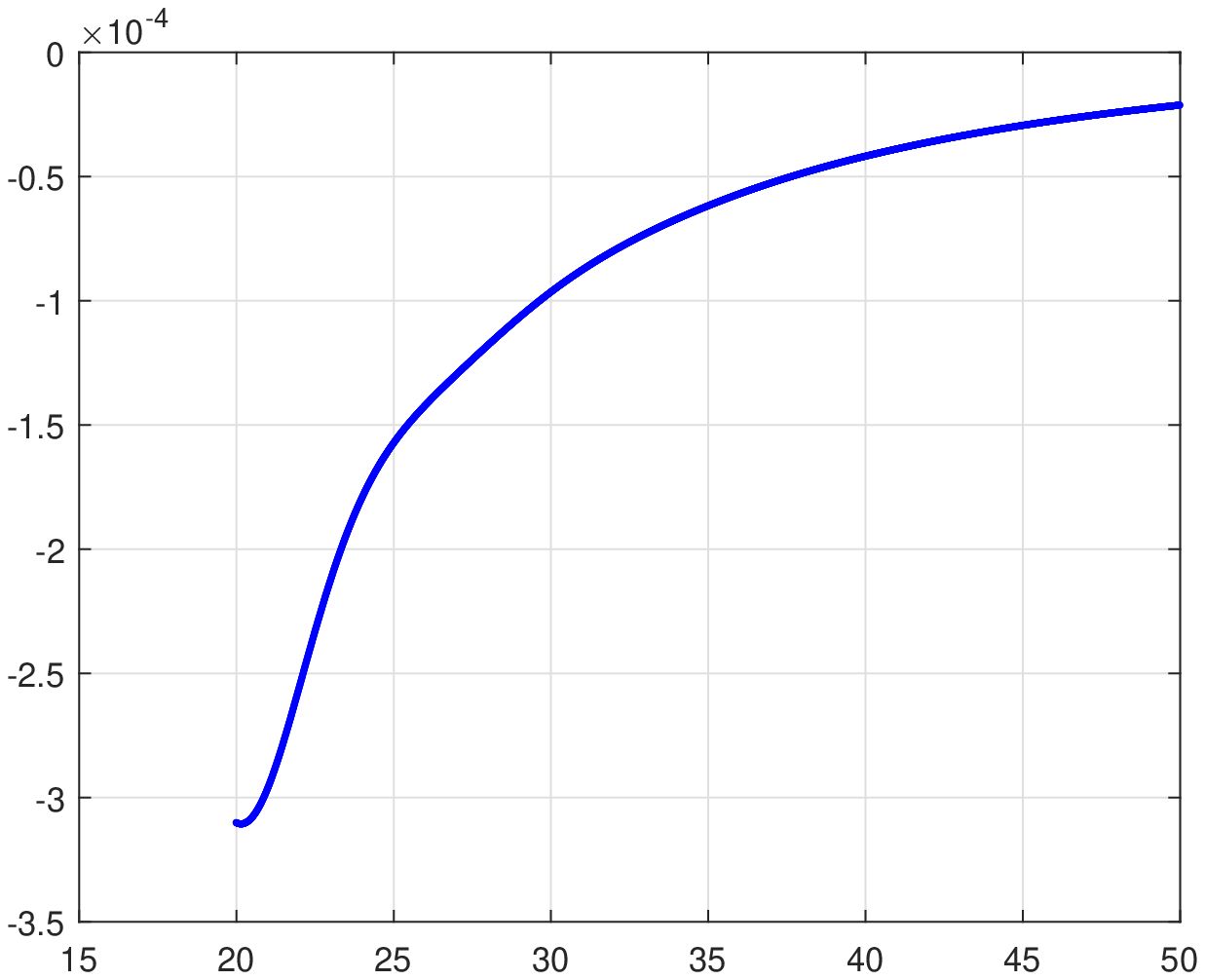}
		\caption{$f(x)=xe^{-x}$. Left: full solution, Right: zoom -- matches \cite{Podlubny}.}\label{ExampleFig2}
	\end{center}
\end{figure}

\begin{figure}[H]
	\begin{center}
		\includegraphics[width=7.5cm]{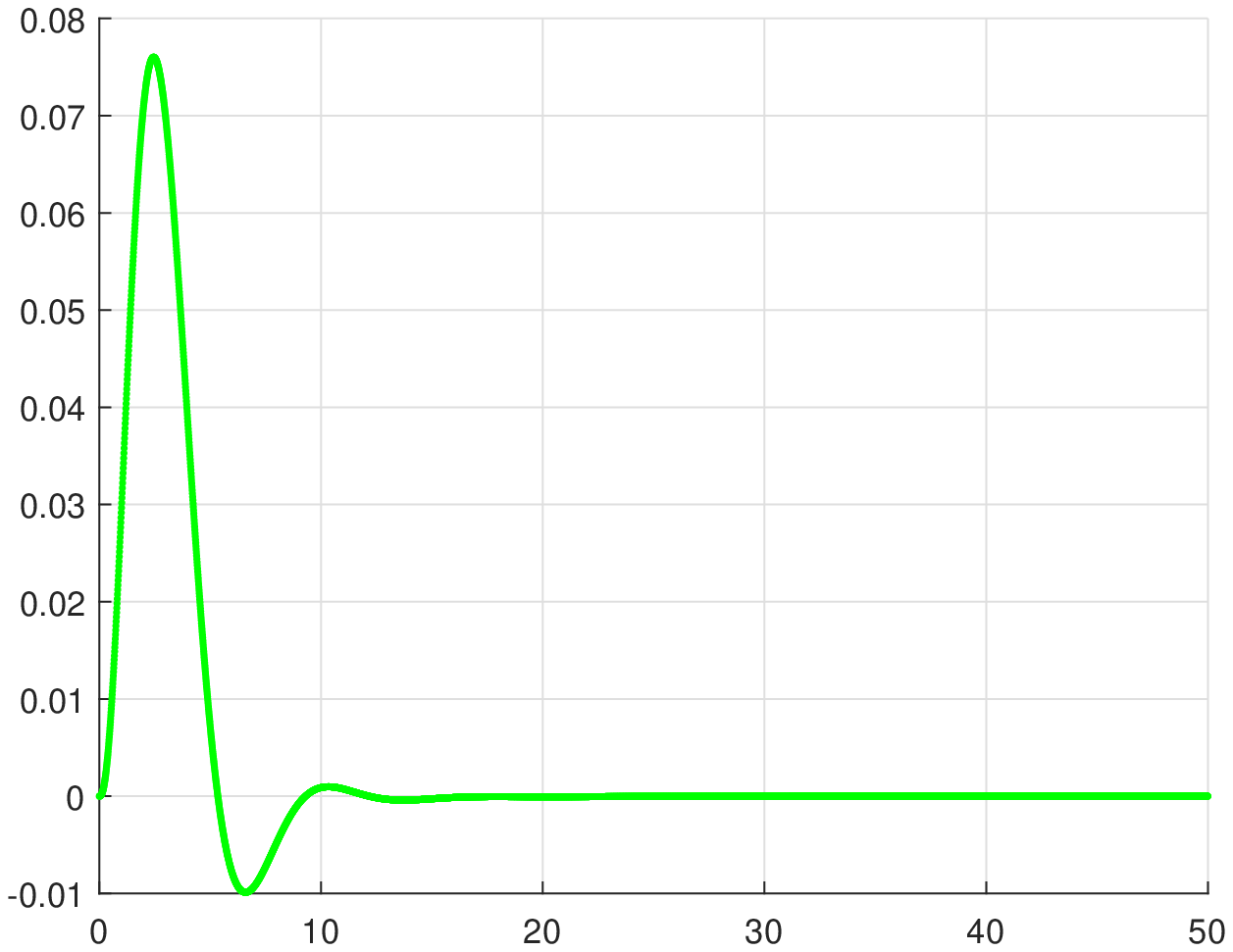}
		\includegraphics[width=7.5cm]{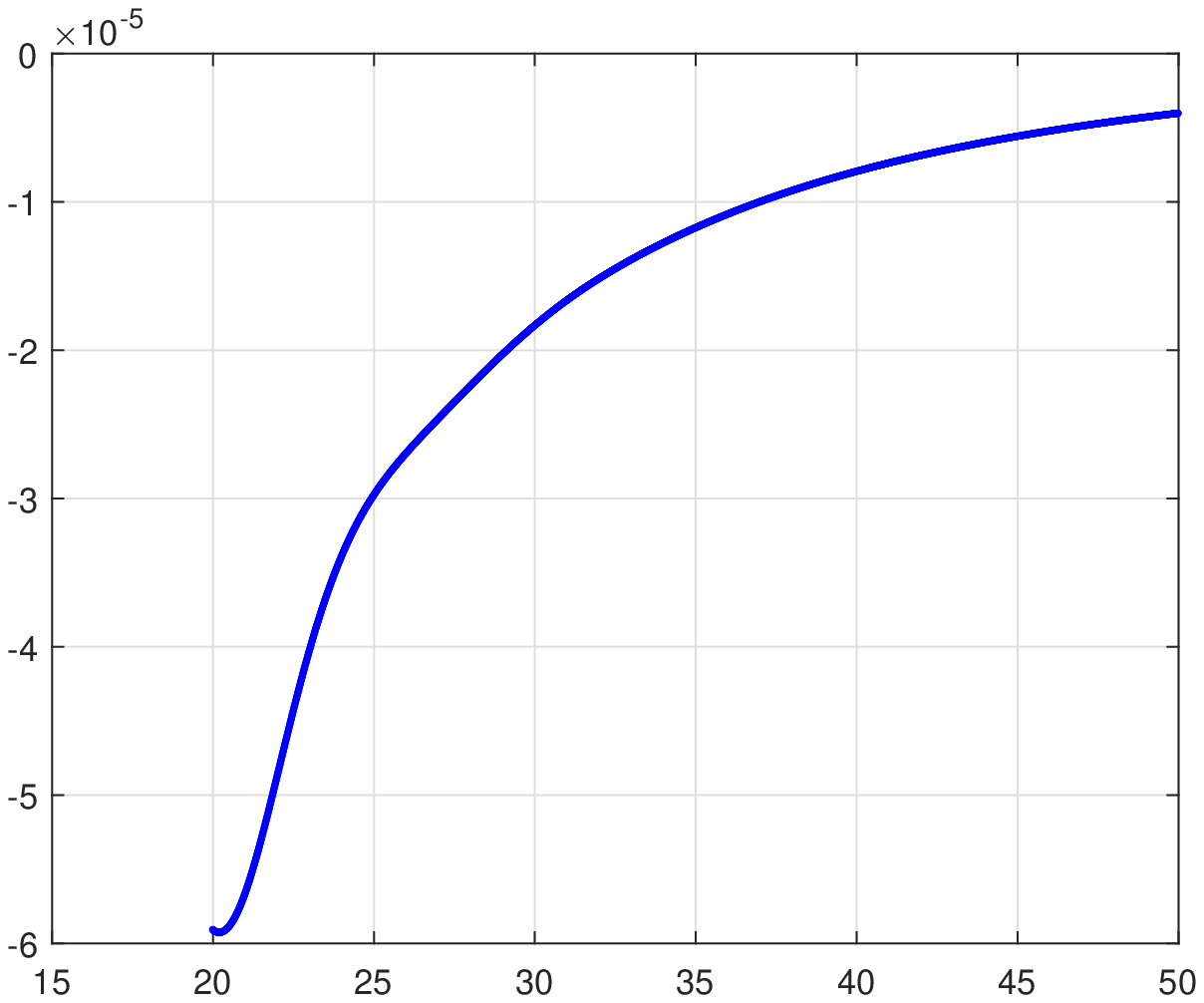}
		\caption{$f(x)=e^{-x}\sin(0.2x)$. Left: full solution, Right: zoom -- matches \cite{Podlubny}.}\label{ExampleFig4}
	\end{center}
\end{figure}

\textbf{Examples 1-4}\\ 
To demonstrate the usage of the substitution method, we consider the initial value problem
$D^{3/2}y(x) + y(x)=f(x)$ with zero initial conditions $y(0)=y'(0)=0$. We consider several examples 
taken from \cite{Podlubny}, page 242, where these problems were solved using Gr\"{u}nwald-Letnikov definition of fractional derivative. The results must be the same because for zero initial data Caputo derivatives coincide with Riemann-Liouville and Gr\"{u}nwald-Letnikov definitions. Our substitution method produces the same results for three options of the right-hand side $f(x)=1$ (Figure \ref{ExampleFig1}), 
$f(x)=xe^{-x}$ (Figure \ref{ExampleFig2}), and $f(x)=e^{-x}\sin(0.2x)$ (Figure \ref{ExampleFig4}).

\begin{figure}[H]
	\begin{center}
		\includegraphics[width=7.5cm]{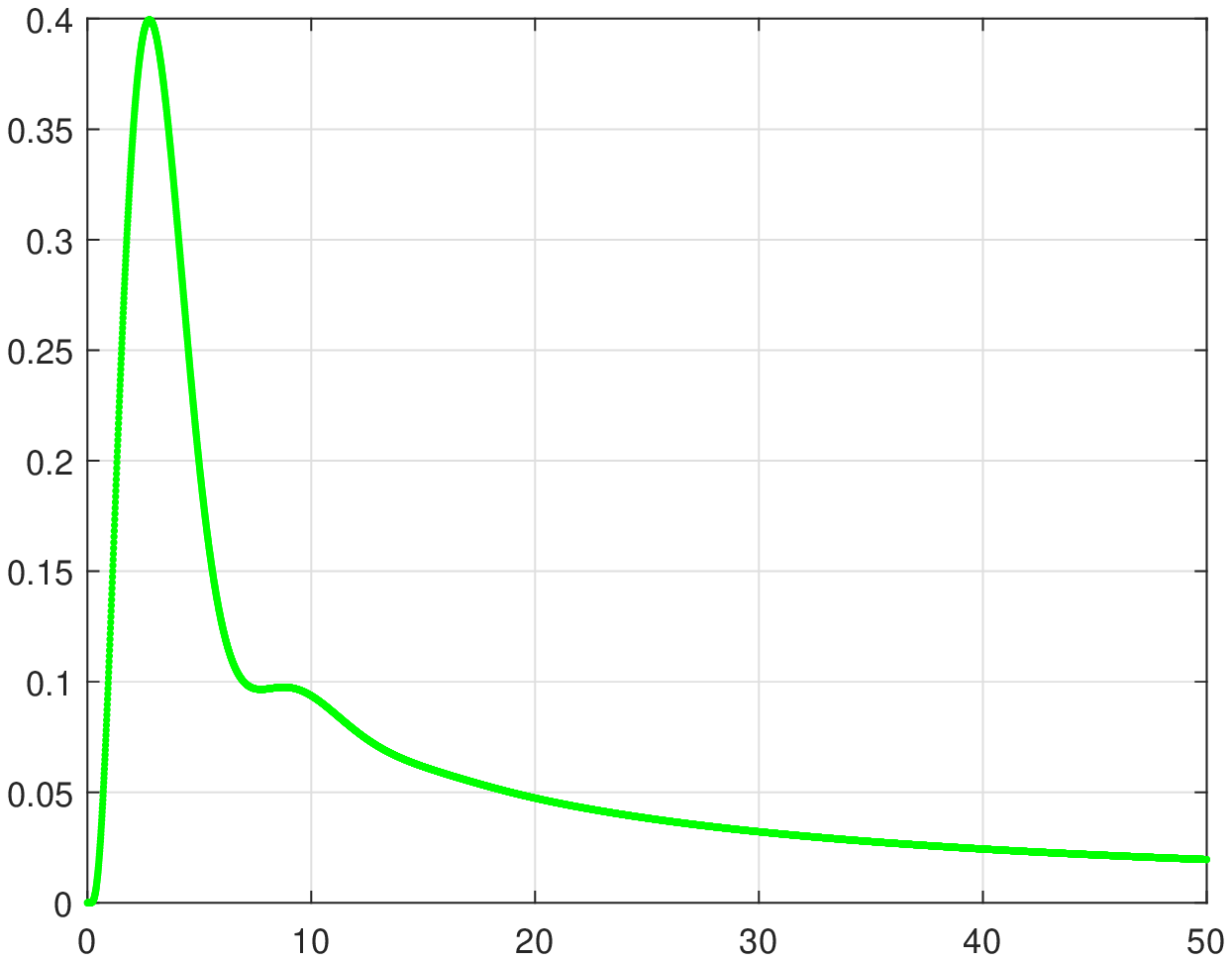}
		\includegraphics[width=7.5cm]{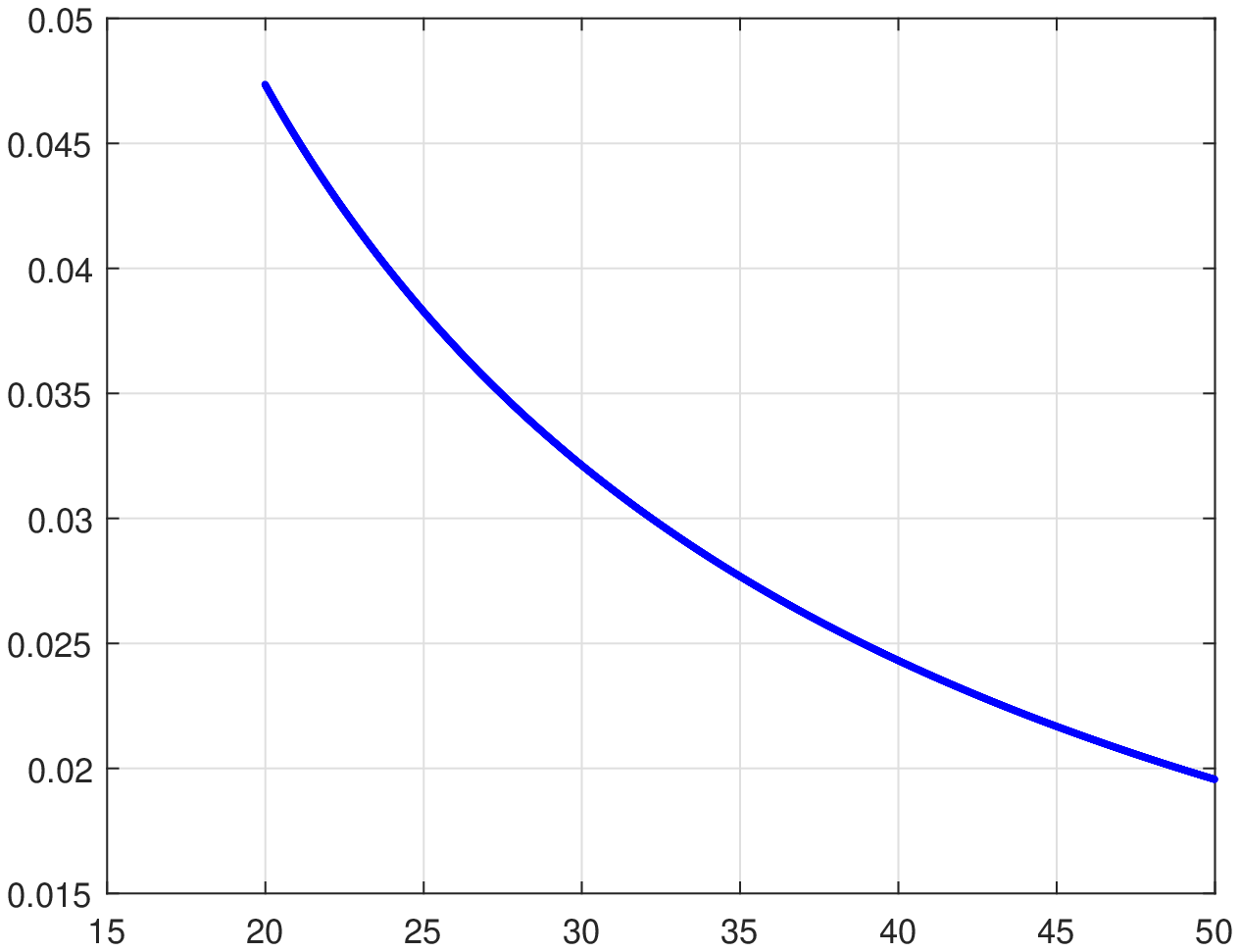}
		\caption{$f(x)=x^{-1}e^{-1/x}$. Left: full solution, Right: zoom.}\label{ExampleFig3}
	\end{center}
\end{figure}

However, in the case $f(x)=x^{-1}e^{-1/x}$ (Figure \ref{ExampleFig3}), our result differs from that in 
\cite{Podlubny} but coincides with another computational method, the method of integration by parts, and, thus, in view of \cite{2020-1Du-Sl}, should be treated as reliable.  

\noindent \textbf{Example 5}\\
Let's consider the extension of the generalized fractional Bessel equation
\cite{2020-4Du-Sl}:
\begin{equation}\label{Example1}
1.5x^{1.5}D^{1.5}u(x)-1.2x^{1.9}D^{1.1}u(x) + 3xD^{0.5}u(x)+(x^2-\nu^2)u(x)=0.
\end{equation}
It can be solved through a series (see \cite{2020-4Du-Sl}): $u(x)=\ds\sum_{n=0}^{\infty}c_n x^{\gamma+sn}$, where $s=0.1$ and for $\nu=2$ we find $\gamma=2.1995$, and for $\nu=3.5$ we have $\gamma=4.3181$ with 4 digits of precision.  

\begin{figure}[H]
	\begin{center}
		\includegraphics[width=7.5cm]{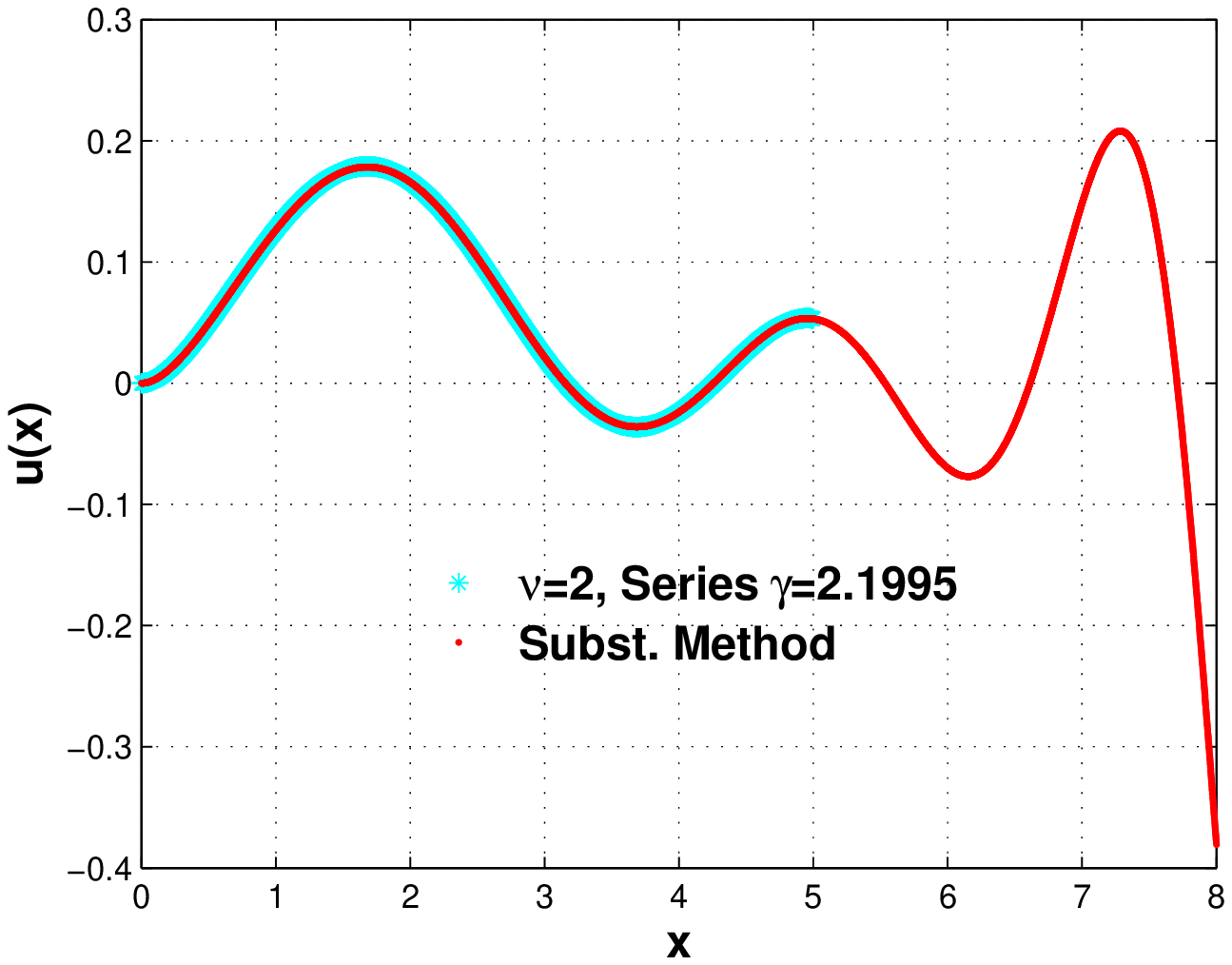}
		\includegraphics[width=7.5cm]{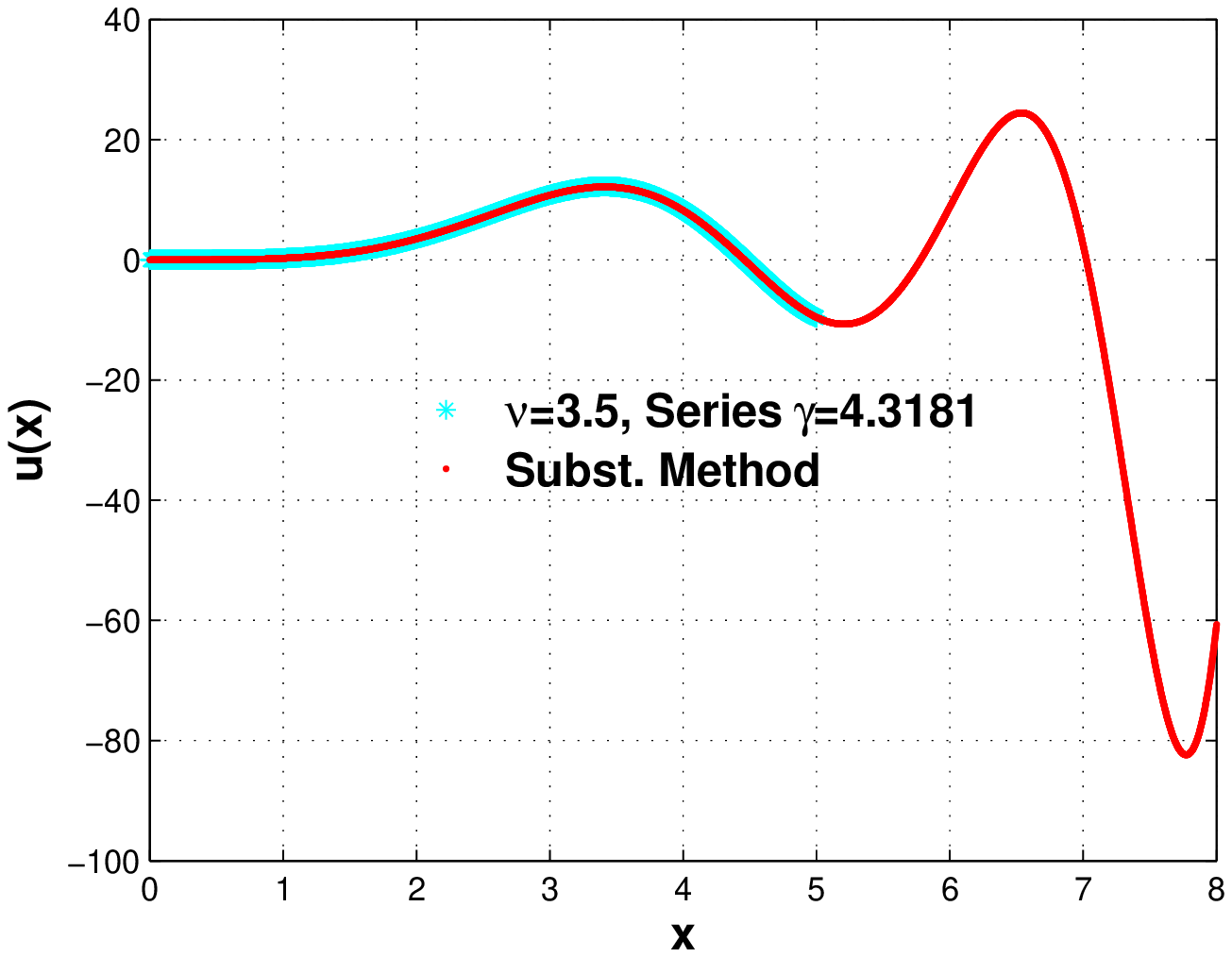}
		\caption{Light Blue - series solution, red - solution generated by the substitution method. On [0,5] graphs of both methods match perfectly.}\label{SolEx1Subst}	
	\end{center}
\end{figure}

Let's point out that for zero initial conditions the substitution method finds only the trivial solution. Therefore we have to disturb the problem by assigning a small value for the derivative. As we know, it affects only the constant in front of function $u(x)$.  After 'calibration' of the first derivative we can generate the solution on a larger interval.  These solutions are shown in Figure \ref{SolEx1Subst}. As we see, the computational solution found by substitution coincides with the exact fractional series solution. 

\section{Conclusions}
\begin{itemize}
\item We suggest the method of substitution for Caputo fractional derivative and its justification.
It provides another way to look at the Caputo derivative and naturally leads to its numerical implementation.
\item We introduce the numerical approximation of the Caputo derivative and derive the method to solve fractional differential equations based on substitution.
\item We present discretization of linear fractional differential equations of any order with utilization of the substitution method.
\item As a part of the analysis of the numerical implementation of the substitution method, we provide a proof of the sufficient condition for well-conditioned problem for linear fractional differential equation.
\item We present examples of computations and demonstrate their fit to the exact solutions and the solutions found by other numerical methods. 
\end{itemize}

\nocite{*}

\end{document}